\documentclass[12pt]{article}
\usepackage{amsmath}
\usepackage{amssymb}
\usepackage[cp1251]{inputenc}
\usepackage[russian]{babel}
\newcommand{\il}[2]{\int\limits_{#1}^{#2}}
\newcommand{\ilp}[1]{\int\limits_{#1}^{+\infty}}

\newcommand{\ph}{\phantom{a}}
\newcommand{\phh}{\phantom{aaa}}

\newcommand{\sist}[2]{\left\{
\begin{array}{l}
{#1}\\
\ph\\
{#2}
\end{array}
\right.}

\begin{document}

MSC 34D20

\vskip 20pt

\centerline{\bf On the stability of systems of two linear first-order}
\centerline{\bf ordinary differential equations}

\vskip 10 pt

\centerline{\bf G. A. Grigorian}

\vskip 10 pt

\centerline{0019 Armenia c. Yerevan, str. M. Bagramian 24/5}
\centerline{Institute of Mathematics NAS of Armenia}
\centerline{E - mail: mathphys2@instmath.sci.am, \ph phone: 098 62 03 05, \ph 010 35 48 61}

\vskip 20 pt

\noindent
Abstract. The Riccati equation method is used to establish some new stability  criteria for systems of two linear first-order ordinary differential equations. It is shown that two of these criteria in the two dimensional case imply
the   Routh - Hurwitz's  criterion.

\vskip 20 pt


\vskip 20 pt
{\bf  1. Introduction}.
Let $a(t), \ph b(t), \ph c(t)$ and $d(t)$ be complex-valued continuous functions on $[t_0,+\infty)$. Consider the linear system
$$
\sist{\phi' = a(t) \phi + b(t) \psi,} {\psi' = c(t) \phi + d(t) \psi, \ph t \ge t_0.} \eqno (1.1)
$$

{\bf Definition 1.1}. {\it A normal linear system of ordinary differential equations (in particular the system (1.1)) is called Lyapunov (asymptotically) stable if its all solutions are bounded (vanish at $+\infty$).}

Study of the stability behavior of the system (1.1), in general, of linear systems of ordinary differential equations is an important problem of Qualitative theory of differential equations, and many works are devoted to it (see [1] and cited works therein, [2 - 4]). The fundamental thorem of R. Bellman (see [5], pp. 168, 169) reduces the study of boundedness conditions of solutions of a wide class of nonlinear systems of ordinary differential equations  to the study of stability conditions of linear systems of ordinary differential equations. There exists various methods of detection of stable and (or) unstable linear systems of ordinary differential equations. Among them notice the Lyapunov, Bogdanov, and Wazevski's methods, the method involving estimates of solutions in the Lozinski's logarithmic norms, and the freezing method (see [Adr], pp. 40 -98). These and other methods (see e. g.; [6 - 10]) permit to carry out wide classes of stable and (or) unstable linear systems. However they cannot completely describe all stable and unstable linear systems of ordinary differential equations in terms of their coefficients.

In this paper on the basis of results of works [11] and [12] by the use of Riccati equation method   new stability criteria for the system (1.1)   are obtained. It is shown that in the two dimensional case of linear systems the  Routh - Hurwitz's   stability criterion is a consequence of the obtained results.

{\bf 2. Auxiliary propositions.} Let $p(t)$ and $q(t)$ be complex-valued continuous functions on $[t_0, + \infty)$. Consider the second order linear ordinary differential equation
$$
\phi'' + p(t) \phi' + q(t) \phi = 0, \phh t\ge t_0. \eqno (2.1)
$$
The substitution $\phi'= \psi$ in this equation reduces it into the linear system
$$
\sist{\phi' = \psi,} {\psi' = - q(t) \phi - p(t) \psi, \ph t\ge t_0.} \eqno (2.2)
$$

{\bf Definition 2.1.} {\it Eq. (2.1) is called Lyapunov (asymptotically)  stable if the corresponding system (2.2) is Lyapunov (asymptotically) stable.}

{\bf Remark 2.1.} {\it It follows from Definition 2.1 that Eq. (2.1) is Lyapunov (asymptotically) stable if and only if its all solutions $\phi(t)$ with $\phi'(t)$ are bounded (vanish at $+ \infty$).}

Set: $G(t) \equiv q(t) - \frac{p'(t)}{2} - \frac{p^2(t)}{4}, \ph \mathcal{L}_0(t) \equiv \frac{1}{\sqrt[4]{G(t)}}\il{t_0}{t}\frac{|(\sqrt{G(\tau)})'|}{\sqrt[4]{G(\tau)}}d \tau, \ph t\ge t_0.$ Hereafter we will assume that $p(t)$ and $G(t)$ are continuously differentiable on $[t_0,+\infty)$, and $G(t) \neq 0, \linebreak t\ge ~t_0.$

{\bf Theorem 2.1.} {\it Let the following conditions be satisfied.

$G(t) > 0, \ph t\ge t_0, \lim\limits_{t \to + \infty}\frac{G'(t)}{G^{3/2}(t)} = \alpha, \ph |\alpha| < 4, \ph \mathcal{L}_0(t)$ and $ Var_{t_0}^t \frac{G'(t)}{G^{3/2}(t)}$ are bounded. Then all solutions of Eq. (2.1) are bounded (vanish at $+\infty$) if and only if
$$
\inf\limits_{t\ge t_0}\biggl\{\il{t_0}{t} \mathfrak{Re} \hskip 2pt p(\tau) d \tau + \frac{1}{2} \ln G(t)\biggr\} > - \infty \phh \biggr( \lim\limits_{t \to + \infty} \biggl\{\il{t_0}{t} \mathfrak{Re} \hskip 2pt p(\tau) d \tau + \frac{1}{2} \ln G(t)\biggr\}  = + \infty\biggr).
$$}

See the proof in [11].

{\bf Theorem 2.2.} {\it Let the conditions of Theorem 2.1 be satisfied. Then Eq (2.1) is Lyapunov (asymptotically) stable if and only if
$$
\sist{\inf\limits_{t\ge t_0}\biggl\{\il{t_0}{t} \mathfrak{Re} \hskip 2pt p(\tau) d \tau  - 2\ln(1 + |p(t)|) + \frac{1}{2} \ln G(t)\biggr\} > - \infty }{\inf\limits_{t\ge t_0}\biggl\{\il{t_0}{t} \mathfrak{Re} \hskip 2pt p(\tau) d \tau - \frac{1}{2} \ln G(t)\biggr\} > - \infty }
$$
$$
\left(\sist{\lim\limits_{t \to + \infty}\biggl\{\il{t_0}{t} \mathfrak{Re} \hskip 2pt p(\tau) d \tau  - 2\ln(1 + |p(t)|) + \frac{1}{2} \ln G(t)\biggr\} = + \infty }{\lim\limits_{t\to + \infty}\biggl\{\il{t_0}{t} \mathfrak{Re} \hskip 2pt p(\tau) d \tau - \frac{1}{2} \ln G(t)\biggr\} = + \infty }\right)
$$
}

See the proof in [11].

For any positive and continuously differentiable  on $[t_0,+\infty)$ function $x(t)$ denote
$$
R_x(t_1;t) \equiv \frac{1 + \sqrt{x(t_0)}(t_1 - t_0)}{1 + \sqrt{x(t_0) (t - t_0)}} \exp\biggl\{-\il{t_1}{t}\sqrt{x(s)} d s\biggr\} \sup\limits_{\xi \in [t_0,t_1]}\frac{|(\sqrt{x(\xi)})'|}{\sqrt{x(\xi)}} + \sup\limits_{\xi \in [t_1,t]}\frac{|(\sqrt{x(\xi)})'|}{\sqrt{x(\xi)}},
$$
$t_0 \le t_1 \le t.$ Set $\rho_x(t) \equiv \inf\limits_{t_1 \in [t_0,t]} R_x(t_1;t), \ph t \ge t_0.$

{\bf Theorem 2.3.} {\it Let the conditions

A) $G(t) < 0, \ph t\ge t_0, \ph p(t)$ and $G(t0$ are continuously differentiable,

\noindent
and one of the following groups of conditions

B) $G(t)$ is non increasing; for some $\varepsilon > 0$ the function $\frac{G'(t)}{|G(t)|^{3/2 -\varepsilon}}$ is bounded;

C) $ - G(t) \ge \varepsilon > 0$; the function $\frac{G'(t)}{G(t)}$ is bounded and
$\ilp{t_0}\rho_{|G|}(\tau)\frac{|G'(\tau)|}{|G(\tau)|^{3/2}} d \tau < + \infty$

\noindent
be satisfied. Then all solutions of Eq. (2.1) are bounded (vanish at $+\infty$) if and only if
$$
\inf\limits_{t \ge t_0}\biggl[\il{t_0}{t}\biggl(\mathfrak{Re}\hskip 2pt p(\tau) - 2 \sqrt{|G(\tau)|}\biggr) + \frac{1}{2} \ln |G(t)|\biggr] > - \infty
$$
$$
\biggl(\lim\limits_{t \to + \infty}\biggl[\il{t_0}{t}\biggl(\mathfrak{Re}\hskip 2pt p(\tau) - 2 \sqrt{|G(\tau)|}\biggr) + \frac{1}{2} \ln |G(t)|\biggr] = + \infty\biggr).
$$
}
 See the proof in [12].

{\bf Theorem 2.4.} {\it Let the conditions A) and the group of conditions C) or the group of conditions

D) $G(t)$ is non increasing, $\frac{G'(t)}{G(t)}$ is bounded

\noindent
be satisfied. Then Eq. (2.1) is Lyapunov (asymptotically) stable if and only if
$$
\inf\limits_{t\ge t_0}\biggl[\il{t_0}{t}\biggl( \mathfrak{Re} \hskip 2 pt p(\tau) - 2\sqrt{|G(\tau)|}\biggr) d \tau + \frac{1}{2} \ln |G(t)| - 2 \ln (1 + |p(t) - 2\sqrt{|G(t)|}\hskip 1pt|)\biggr] > - \infty
$$
$$
\biggl(\lim\limits_{t\to +\infty}\biggl[\il{t_0}{t}\biggl( \mathfrak{Re} \hskip 2 pt p(\tau) - 2\sqrt{|G(\tau)|}\biggr) d \tau + \frac{1}{2} \ln |G(t)| - 2 \ln (1 + |p(t) - 2\sqrt{|G(t)|}\hskip 1pt|)\biggr] = + \infty\biggr).
$$
}
 See the proof in [12].

{\bf Corollary 2.1.} {\it Assume $-G(t) \ge \varepsilon > 0, \ph t\ge t_0; \ph \frac{|G'(t)}{|G(t)}| \le \frac{M}{(1 +t - t_0)^\alpha}, \ph t\ge t_0, \ph M >~ 0, \linebreak \alpha > 0, \ph \ilp{t_0}\frac{d \tau}{\sqrt{|G(\tau)|}(1 + \tau - t_0)^{2\alpha}} < + \infty$ and let the conditions A) be satisfied. Then the following statements are valid.

$A_1)$ All solutions of Eq. (2.1) are bounded (vanish at $+\infty$) if and only if
$$
\inf\limits_{t \ge t_0}\biggl[\il{t_0}{t}\biggl(\mathfrak{Re}\hskip 2pt p(\tau) - 2 \sqrt{|G(\tau)|}\biggr) + \frac{1}{2} \ln |G(t)|\biggr] > - \infty
$$
$$
\biggl(\lim\limits_{t \to + \infty}\biggl[\il{t_0}{t}\biggl(\mathfrak{Re}\hskip 2pt p(\tau) - 2 \sqrt{|G(\tau)|}\biggr) + \frac{1}{2} \ln |G(t)|\biggr] = + \infty\biggr);
$$

$B_1)$ Eq. (2.1) is Lyapunov (asymptotically) stable if and only if
$$
\inf\limits_{t\ge t_0}\biggl[\il{t_0}{t}\biggl( \mathfrak{Re} \hskip 2 pt p(\tau) - 2\sqrt{|G(\tau)|}\biggr) d \tau + \frac{1}{2} \ln |G(t)| - 2 \ln (1 + |p(t) - 2\sqrt{|G(t)|}\hskip 1pt|)\biggr] > - \infty
$$
$$
\biggl(\lim\limits_{t\to +\infty}\biggl[\il{t_0}{t}\biggl( \mathfrak{Re} \hskip 2 pt p(\tau) - 2\sqrt{|G(\tau)|}\biggr) d \tau + \frac{1}{2} \ln |G(t)| - 2 \ln (1 + |p(t) - 2\sqrt{|G(t)|}\hskip 1pt|)\biggr] = + \infty\biggr).
$$
}
 See the proof in [12].

Consider the Riccati equations
$$
y' + b(t) y^2 + A(t) y - c(t) = 0, \phh t\ge t_0, \eqno (2.3)
$$
$$
z' + c(t) z^2 - A(t) z - a(t) = 0, \phh t\ge t_0, \eqno (2.4)
$$
where $A(t) \equiv a(t) - d(t), \ph t\ge t_0$. It is not difficult to verify that the solutions $y(t) \ph (z(t))$ of Eq. (2.3) (Eq. (2.4)), existing on an interval $[t_1,t_2) \ph (t_0 \le t_1 < t_2 \le + \infty)$ are connected with solutions $(\phi(t), \psi(t))$ of the system (1.1) by relations (see e.g.; [2])
$$
\phi(t)= \phi(t_1) \exp\biggl\{\il{t_1}{t}\bigl[b(\tau) y(\tau) + a(\tau)\bigr] d \tau\biggr\}, \ph \phi(t_1) \ne 0, \ph \psi(t) = y(t) \phi(t), \ph t\in [t_1, t_2) \eqno (2.5)
$$
$$
\biggl(\psi(t)= \psi(t_1) \exp\biggl\{\il{t_1}{t}\bigl[c(\tau) z(\tau) + d(\tau)\bigr] d \tau\biggr\}, \ph \psi(t_1) \ne 0, \ph \phi(t) = z(t) \psi(t), \biggr) \eqno (2.6)
$$
$t\in [t_1,t_2)$. Hereafter we will assume that $a(t), \ph b(t), \ph c(t)$ and $d(t)$ are continuously differentiable on $[t_0,+\infty)$ and $a(t) \ne 0, \ph c(t) \ne 0, \ph t\ge t_0$. Set:
$$
D_1(t) \equiv \frac{a(t) b'(t) - a'(t) b(t)}{b(t)} + a(t) d(t) - b(t) c(t), \phantom{aaaaaaaaaaaaaaaaaaaaaaaaaaaaaaaaaaaaaaaaaaaa}
$$
$$
\phantom{aaaaaaaaaaaaaaaaaaaaaaaaa} D_2(t) \equiv \frac{d(t) c'(t) - d'(t) c(t)}{c(t)} + a(t) d(t) - b(t) c(t), \phh t\ge t_0.
$$
The substitution
$$
u = b(t) y + a(t), \phh t\ge t_0 \eqno (2.7)
$$
in Eq. (2.3) transforms that into the equation
$$
u' + u^2 - \biggl[S(t) + \frac{b'(t)}{b(t)}\biggr] u + D_1(t) = 0, \phh t\ge t_0, \eqno (2.8)
$$
where $S(t) \equiv a(t) + d(t), \ph t\ge t_0$. Analogously the substitution
$$
v = c(t) z + d(t), \phh t\ge t_0 \eqno (2.9)
$$
in Eq. (2.4) transforms that into the equation
$$
v' + v^2 - \biggl[S(t) + \frac{c'(t)}{c(t)}\biggr] v + D_2(t) = 0, \phh t\ge t_0, \eqno (2.10)
$$
Consider the second order linear ordinary differential equations
$$
\phi'' - \biggl[S(t) + \frac{b'(t)}{b(t)}\biggr]\phi' + D_1(t) \phi =0, \phh t \ge t_0, \eqno (2.11)
$$
$$
\psi'' - \biggl[S(t) + \frac{c'(t)}{c(t)}\biggr]\psi' + D_2(t) \psi =0, \phh t \ge t_0. \eqno (2.12)
$$
It is not difficult to verify that the solutions $u(t) \ph (v(t))$ of Eq. (2.8) (Eq (2.10)), existing on $[t_1,t_2)$, are connected with solutions $\phi_0(t), \ph (\psi_0(t))$ of Eq. (2.11) (Rq. (2.12)) by relations
$$
\phi_0(t) = \phi_0(t_1) \exp\biggl\{\il{t_1}{t} u(\tau) d \tau\biggr\}, \ph \phi_0(t_1)\ne 0, \ph t\in [t_1,t_2), \eqno (2.13)
$$
$$
\psi_0(t) = \psi_0(t_1) \exp\biggl\{\il{t_1}{t} v(\tau) d \tau\biggr\}, \ph \psi_0(t_1)\ne 0, \ph t\in [t_1,t_2), \eqno (2.14)
$$
On the other hand by (2.5) - (2.7) and (2.9) the same solutions $u(t)$ and $v(t)$ are connected with solutions $(\phi(t), \psi(t))$ of the system (1.1) by relations
$$
\phi(t) = \phi(t_1) \exp\biggl\{\il{t_1}{t} u(\tau) d \tau\biggr\}, \phh \psi(t) = \psi(t_1) \exp\biggl\{\il{t_1}{t} v(\tau) d \tau\biggr\}, \phh t\in [t_1,t_2), \eqno (2.15)
$$
$
\phi(t_1) \ne 0, \ph \psi(t_1) \ne 0, \ph \frac{u(t_1) - a(t_1)}{b(t_1)} \hskip 2pt \frac{v(t_1) - d(t_1)}{c(t_1)} = 1.
$
By (2.5) - (2.7) and (2.9) the last equality is equivalent to the following one
$$
\biggl[\frac{\phi'(t_1)}{\phi(t_1)} - a(t_1)\biggr] \biggl[\frac{\psi'(t_1)}{\psi(t_1)} - d(t_1)\biggr]  = b(t_1) c(t_1). \eqno (2.16)
$$
By the uniqueness theorem from (2.13) - (2.16) we immediately get

{\bf Lemma 2.1.} {\it Let $\phi_0(t)$ and $\psi_0(t)$ be solutions of Eq. (2.11) and (2.12) respectively such that $\phi_0(t) \ne 0, \ph \psi_0(t) \ne 0, \ph t\in [t_1,t_2), \ph \biggl[\frac{\phi_0'(t_1)}{\phi_0(t_1)} - a(t_1)\biggr] \biggl[\frac{\psi_0'(t_1)}{\psi_0(t_1)} - d(t_1)\biggr]  = b(t_1) c(t_1).$ Then $(\phi_0(t), \psi_0(t))$ is a solution of the system (1.1) on $[t_1,t_2)$.}

$\phantom{aaaaaaaaaaaaaaaaaaaaaaaaaaaaaaaaaaaaaaaaaaaaaaaaaaaaaaaaaaaaaaaaaaaa} \Box$

Hereafter we will assume that $S(t) + \frac{b'(t)}{b(t)}$ and $S(t) + \frac{c'(t)}{c(t)}$ are continuously differentiable on $[t_0, + \infty)$.  Set:
$$
G_1(t) \equiv D_1(t) + \frac{1}{2} \biggl[S(t) + \frac{b'(t)}{b(t)}\biggr]' - \frac{1}{4} \biggl[S(t) + \frac{b'(t)}{b(t)}\biggr]^2, \phh t\ge t_0,
$$
$$
G_2(t) \equiv D_2(t) + \frac{1}{2} \biggl[S(t) + \frac{c'(t)}{c(t)}\biggr]' - \frac{1}{4} \biggl[S(t) + \frac{c'(t)}{c(t)}\biggr]^2, \phh t\ge t_0.
$$

{\bf Lemma 2.2.} {\it Assume $Im \hskip 2pt G_1(t) \equiv 0 \ph (Im \hskip 2pt G_2(t) \equiv 0), \ph t\ge t_0$, and $Im \hskip 2pt \biggl[\lambda - \frac{1}{2}\biggl(S(t_0) + \frac{b'(t_0)}{b(t_0)}\biggr)\biggr] \ne 0 \ph \biggl(Im \hskip 2pt \biggl[\lambda - \frac{1}{2}\biggl(S(t_0) + \frac{c'(t_0)}{c(t_0)}\biggr)\biggr] \ne 0\biggr) $ for some complex $\lambda$. Then Eq. (2.8) (Eq. (2.10)) has a solution $u(t) \ph (v(t))$ on $[t_0, +\infty)$ with $u(t_0) = \lambda \ph (v(t_0) = \lambda)$.}

Proof. In Eq. (2.8) substitute
$$
u = w + \frac{1}{2}\biggl(S(t) + \frac{b'(t)}{b(t)}\biggr), \phh t\ge t_0. \eqno (2.17)
$$
We obtain
$$
w' + w^2 + G_1(t) = 0, \phh t\ge t_0.. \eqno (2.18)
$$
Show that this equation has a solution $w(t)$ on $[t_0, +\infty)$ with $w(t_0) = \lambda + \frac{1}{2}\biggl[S(t_0) + \frac{b'(t_0)}{b(t_0)}\biggr].$ Consider the second order linear ordinary differential equation
$$
\chi'' + G_1(t) \chi = 0, \phh t\ge t_0.
$$
Let $\chi_1(t)$ and $\chi_2(t)$ be the solutions of this equation on $[t_0, +\infty)$ with $\chi_k(t_0) =1, \linebreak k=1,2, \ph \chi'_1(t_0) = \lambda_1 - \lambda_2, \ph \chi'_2(t_0) = \lambda_1 + \lambda_2,$ where $\lambda_1 \equiv \mathfrak{Re}\hskip 2pt \biggl[\lambda - \frac{1}{2} \biggl(S(t_0) + \frac{b'(t_0)}{b(t_0)}\biggr)\biggr], \linebreak \lambda_2 \equiv Im \hskip 2pt \biggl[\lambda - \frac{1}{2} \biggl(S(t_0) + \frac{b'(t_0)}{b(t_0)}\biggr)\biggr] \ne 0.$
Since $G_1(t)$ is a real-valued function $\chi_k(t), \ph k=1,2$ are also real-valued ones. Moreover, obviously,  $\chi_k(t), \ph k=1,2$ are linearly independent. Consequently $\chi(t)\equiv \chi_1(t) + i \chi_2(t) \ne 0, \ph t\ge t_0$ and $w(t) \equiv\frac{\chi'(t)}{\chi(t)}$ is a solution of Eq. (2.18) on $[t_0,+\infty)$ with $w(t_0) = \lambda - \frac{1}{2}\biggl(S(t_0) + \frac{b'(t_0)}{b(t_0)}\biggr)$.  Then by $(2.17)$ $u(t) \equiv v(t) + \frac{1}{2}\biggl(S(t) + \frac{b'(t)}{b(t)}\biggr)$ is a solution of Eq. (2.8) on $[t_0,+\infty)$ with $u(t_0) = \lambda.$. Existence of $v(t)$ can be proved by analogy. The lemma is proved.

{\bf Theorem 2.5.} {\it The following statements are valid.

I. The system (1.1) is Lyapunov (asymptotically) stable if and only if all solutions of Eq. (2.11) and Eq. (2.12) are bounded (vanish at $+\infty$).

II. Assume $a(t), \ph b(t)$ and $\frac{1}{b(t)}$ are bounded. Then the system (1.1) is Lyapunov \linebreak (asymptotically) stable if and only if Eq. (2.11) is Lyapunov (asymptotically) stable.}

Proof. Obviously there exist $\lambda_1 \ne \lambda_2$ such that $Im \hskip 2pt \biggl[\lambda_k - \frac{1}{2}\biggl(S(t_0) + \frac{b'(t_0)}{b(t_0)}\biggr)\biggr] \ne 0, \linebreak Im \hskip 2pt \biggl[\frac{b(t_0) c(t_0)}{\lambda_k - a(t_0)} + d (t_0) - \frac{1}{2}\biggl(S(t_0) + \frac{c'(t_0)}{c(t_0)}\biggr)\biggr] \ne 0, \ph k=1,2.$ Let $u_k(t) \ph (v_k(t)), \ph k=1,2$ be solutions of Eq. (3.8) (Eq. (2.10)) with $u_k(t_0) = \lambda_k \ph (v_k(t_0) = \frac{b(t_0)c(t_0)}{\lambda_k - a(t_0)} + d(t_0)), \ph k=1,2.$ Then by Lemma 2.2 $u_k(t) \ph (v_k(t))), \ph k=1,2$ exist on $[t_0,+\infty)$; moreover
$$
[u_k(t_0) - a(t_0)][v_k(t_0)- d(t_0)] = b(t_0) c(t_0), \phh k=1,2. \eqno (2.19)
$$
Set: $\phi_k(t) \equiv \exp\biggl\{\il{t_0}{t} u_k(\tau) d \tau\biggr\}, \ph \psi_k(t) \equiv \exp\biggl\{\il{t_0}{t}v_k(\tau) d\tau\biggr\}, \ph t\ge t_0, \ph k=1,2.$ By (2.13)  (by (2.14)) $\phi_k(t) \ph (\psi_k(t)), \ph k=1,2$ are solutions of Eq. (2.11) (of Eq. (2.12)) on $[t_0,+\infty)$ and by (2.19) we have
$$
\biggl[\frac{\phi_k'(t_0)}{\phi_k(t_0)} - a(t_0)\biggr] \biggl[\frac{\psi_k'(t_0)}{\psi_k(t_0)} - d(t_0)\biggr] = b(t_0) c(t_0), \phh k=1,2.
$$
In virtue of Lemma 2.1 from here it follows that $(\phi_k(t), \psi_k(t)), \ph k=1,2$ are solutions of the system (1.1) on $[t_0,+\infty)$.  Let us prove statement I. Assume all solutions of Eq. (2.11) and (2.12) are bounded (vanish at
$+\infty$). Then the linearly independent solutions $(\phi_k(t), \psi_k(t)), \ph k=1,2$ are bounded (vanish at $+\infty$). Consequently the system (1.1) is Lyapunov (asymptotically) stable. Assume now  the system (1.1) is Lyapunov (asymptoti- \linebreak cally) stable. Then the linearly independent solutions $\phi_k(t) \ph (\psi_k(t)), \ph k=1,2$ of Eq. (2.11) (of Eq. (2.12)) are bounded (vanish at $+\infty$). Therefore all solutions of Eq. (2.11) and  Eq. (2.12) are bounded (vanish at $+\infty$). The statement I is proved.  Prove statement II. Assume Eq. (2.11) is Lyapunov (asymptotically) stable. Then the functions $\phi_k(t), \ph \phi'_k(t), \ph k=~1,2$ are bounded (vanish at $+\infty$).
 Since by (1.1) $\psi_k(t) = - \frac{a(t)}{b(t)} \phi_k(t) + \frac{1}{b(t)} \phi'_k(t), \ph k=1,2$ and $\frac{a(t)}{b(t)}, \ph \frac{1}{b(t)}$ are bounded the functions $\psi_k(t), \ph k=1,2$ are bounded (vanish at $+\infty$) as well. So the linearly independent solutions $(\phi_k(t), \psi_k(t)), \ph k=1,2$ of the system (1.1)  are bounded (vanish at $+\infty$). Therefore the system (1.1) is Lyapunov (asymptotically) stable. Let now  the system (1.1) be Lyapunov (asymptotically) stable. Then the functions $\phi_k(t), \ph \psi_k(t), \ph k=~1,2$ are bounded (vanish at $+\infty)$. Since by (1.1) $\phi'_k(t) = a(t) \phi_k(t) + b(t) \psi_k(t), \ph t\ge t_0$ and the functions $a(t)$ and $b(t)$ are bounded the functions $\phi'_k(t) ,\linebreak k=~1,2$ are also bounded (vanish at $+ \infty$). Thus all solutions $\phi(t)$ of Eq. (2.11) with $\phi'(t)$  are bounded (vanish at $+ \infty$). Therefore Eq. (2.11) is Lyapunov (asymptotically) stable. The theorem is proved.

{\bf Remark 2.2.} {\it From the proof of statement II is seen that the restrictions on $c(t)$ for that statement are not obligatory.}

{\bf 3. Main results.}  In this section we study the stability behavior of the system (1.1) in the following cases

I. $G_1(t) > 0, \phh G_2(t) > 0, \phh t\ge t_0;$

II. $G_1(t) > 0, \phh G_2(t) < 0, \phh t\ge t_0;$

III. $G_1(t) < 0, \phh G_2(t) < 0, \phh t\ge t_0;$

IV. $G_1(t) > 0, \phh t\ge t_0; \phh$ V. $G_2(t) < 0, \phh t\ge t_0.$

\noindent
The case VI. $G_1(t) < 0, \ph G_2(t) > 0, \ph t\ge t_0$ is reducible to the case III by simple transformation $ \phi \to - \phi$.

{\bf Remark 3.1.} {\it It is easy to study the trivial case $G_1(t) = G_2(t) \equiv 0, \ph t\ge t_0$ separately.}

Set:
$$
\mathcal{L}_k(t) \equiv \frac{1}{\sqrt[4]{G_k(t)}}\il{t_0}{t}\frac{|(\sqrt{G_k(\tau)})'|}{\sqrt[4]{G_k(\tau)}} d\tau, \phh k=1,2, \phh t\ge t_0.
$$

{\bf Theorem 3.1.} {\it Let the following conditions be satisfied

1) $G_k(t) > 0, \phh t \ge t_), \phh \lim\limits_{t \to +\infty}\frac{G_k'(t)}{G_k^{3/2}(t)} = \alpha_k, \ph |\alpha_k| < 4, \ph k=1,2;$

2) $\mathcal{L}_k(t)$ and $Var_{t_0}^t\frac{G_k'(t)}{G_k^{3/2}(t)}$ are bounded $k=1,2.$

\noindent
Then the system (1.1) is Lyapunov (asymptotically) stable if and only if
$$
\sist{\sup\limits_{t\ge t_0}\biggl[\il{t_0}{t} \mathfrak{Re} \hskip 2pt S(\tau) d\tau - \ln |b(t)| - \frac{1}{2}\ln G_1(t)\biggr] < + \infty,}{\sup\limits_{t\ge t_0}\biggl[\il{t_0}{t} \mathfrak{Re} \hskip 2pt S(\tau) d\tau - \ln |c(t)| - \frac{1}{2}\ln G_2(t)\biggr] < + \infty.}  \eqno (3.1)
$$
$$
\left(\sist{\lim\limits_{t\to + \infty}\biggl[\il{t_0}{t} \mathfrak{Re} \hskip 2pt S(\tau) d\tau - \ln |b(t)| - \frac{1}{2}\ln G_1(t)\biggr] = - \infty,}{\lim\limits_{t\to + \infty}\biggl[\il{t_0}{t} \mathfrak{Re} \hskip 2pt S(\tau) d\tau - \ln |c(t)| - \frac{1}{2}\ln G_2(t)\biggr] = - \infty.}\right)  \eqno (3.2)
$$
}

Proof. By virtue of Theorem 2.1 from conditions 1) ,2) it follows that the solutions of Eq. (2.11) and (2.12) are bounded (vanish at $+\infty$) if and only if the inequalities (3.1) (the equalities (3.2)) are satisfied. Then by statement I of Theorem 2.5 the system (1.1) is Lyapunov (asymptotically) stable if and only if the inequalities (3.1) (the equalities (3.2)) are fulfilled. The theorem is proved.

\pagebreak
{\bf Theorem 3.2}. {\it Let the following conditions be satisfied

3) $G_1(t) > 0, \ph t\ge t_0, \phh \lim\limits_{t \to +\infty}\frac{G_1'(t)}{G_1^{3/2}(t)} = \alpha_1, \ph |\alpha_1| < 4;$

4) $\mathcal{L}_1(t)$ and $Var_{t_0}^t\frac{G_1'(t)}{G_1^{3/2}(t)}$ are bounded;

5) $G_2(t) < 0, \ph t\ge t_0,$ and is non increasing, $\frac{G_2'(t)}{|G_2(t)|^{3/2 - \varepsilon}}$ is bounded for some $\varepsilon >0$, or

$5_1) \ph |G_2(t)| \ge \varepsilon > 0, \ph \frac{G_2'(t)}{G_2(t)}$ is bounded and $\ilp{t_0}\rho_{|G_2|}(\tau)\frac{|G_2'(\tau)|}{|G_2(\tau)|^{3/2}} d \tau < + \infty$.

\noindent
Then the system (1.1) is Lyapunov (asymptotically) stable if and only if
$$
\sist{\sup\limits_{t\ge t_0}\biggl[\il{t_0}{t} \mathfrak{Re} \hskip 2pt S(\tau) d\tau - \ln |b(t)| - \frac{1}{2}\ln G_1(t)\biggr] < + \infty,}{\sup\limits_{t\ge t_0}\biggl[\il{t_0}{t} \biggl(\mathfrak{Re} \hskip 2pt S(\tau) + \sqrt{|G_2(\tau)|}\biggr) d\tau - \ln |c(t)| - \frac{1}{2}\ln |G_2(t)|\biggr] < + \infty.}  \eqno (3.3)
$$
$$
\left(\sist{\lim\limits_{t\to + \infty}\biggl[\il{t_0}{t} \mathfrak{Re} \hskip 2pt S(\tau) d\tau - \ln |b(t)| - \frac{1}{2}\ln G_1(t)\biggr] = - \infty,}{\lim\limits_{t\to + \infty}\biggl[\il{t_0}{t}\biggl( \mathfrak{Re} \hskip 2pt S(\tau)  + \sqrt{|G_2(\tau)|}\biggr) d\tau - \ln |c(t)| - \frac{1}{2}\ln |G_2(t)|\biggr] = - \infty.}\right)  \eqno (3.4)
$$
}

Proof. By Theorem 2.1 from conditions 3), 4) it follows that the solutions of Eq. (2.11) are bounded (vanish at $+\infty$) if and only if the first of the inequalities (3.3) (the first of the equalities (3.4)) is satisfied. By Theorem 2.3 from conditions 5) or $5_1$) it follows that the solutions of Eq. (2.12) are bounded (vanish at $+\infty$) if and only if the second of the inequalities(3.3) ( of the equalities (3.4)) is satisfied. Then by Theorem 2.5 the system (1.1) is Lyapunov (asymptotically) stable if and only if the inequalities (3.3) (the equalities (3.4)) are satisfied. The theorem is proved.

By analogy can be proved

{\bf Theorem 3.3.} {\it Let the following conditions be satisfied

6) $G_k(t) < 0, \ph t\ge t_0, \ph G_k(t)$ is non increasing $k=1,2$

7) $\frac{G_k'(t)}{|G_k(t)|^{3/2 - \varepsilon}}$ is bounded for some $\varepsilon > 0, \ph k=1,2$ or

$7_1$) $|G_k(t)| \ge \varepsilon > 0, \ph t\ge t_0, \ph \frac{G_k'(t)}{G_k(t)}$ is bounded and $$\ilp{t_0}\rho_{|G_k|}(\tau)\frac{|G_k'(\tau)|}{|G_k(\tau)|^{3/2}} d \tau < + \infty, \ph k=1,2$$.
Then the system (1.1) is Lyapunov (asymptotically) stable if and only if
$$
\sist{\sup\limits_{t\ge t_0}\biggl[\il{t_0}{t} \biggl(\mathfrak{Re} \hskip 2pt S(\tau) d\tau  + \sqrt{|G_1(\tau)|}\biggr)   - \ln |b(t)| - \frac{1}{2}\ln |G_1(t)|\biggr] < + \infty,}{\sup\limits_{t\ge t_0}\biggl[\il{t_0}{t} \biggl(\mathfrak{Re} \hskip 2pt S(\tau) + \sqrt{|G_2(\tau)|}\biggr) d\tau - \ln |c(t)| - \frac{1}{2}\ln |G_2(t)|\biggr] < + \infty.}  \
$$
$$
\left(\sist{\lim\limits_{t\to + \infty}\biggl[\il{t_0}{t} \biggl(\mathfrak{Re} \hskip 2pt S(\tau) + \sqrt{|G_1(\tau)|}\biggr) d\tau - \ln |b(t)| - \frac{1}{2}\ln |G_1(t)|\biggr] = - \infty,}{\lim\limits_{t\to + \infty}\biggl[\il{t_0}{t}\biggl( \mathfrak{Re} \hskip 2pt S(\tau)  + \sqrt{|G_2(\tau)|}\biggr) d\tau - \ln |c(t)| - \frac{1}{2}\ln |G_2(t)|\biggr] = - \infty.}\right)
$$
}
$\phantom{aaaaaaaaaaaaaaaaaaaaaaaaaaaaaaaaaaaaaaaaaaaaaaaaaaaaaaaaaaaaaaaaaaaaa} \Box$

{\bf Theorem 3.4.} {\it Let the following conditions be satisfied

8) $a(t), \ph b(t)$ and $\frac{1}{b(t)}$ are bounded;

9) $G_1(t) > 0, \ph t\ge t_0, \ph \lim\limits_{t\to +\infty}\frac{G_1'(t)}{G_1^{3/2}(t)} = \alpha_1, \ph |\alpha_1| < 4$;

10)  $\mathcal{L}_1(t)$ and $Var _{t_0}^i \frac{G_1'(t)}{G_1^{3/2}(t)}$ are bounded.

\noindent
Then the system (1.1) is Lyapunov (asymptotically) stable if and only if
$$
\sist{\sup\limits_{t\ge t_0}\biggl[\il{t_0}{t} \mathfrak{Re} \hskip 2pt S(\tau) d\tau + \ln |b(t)| - \frac{1}{2}\ln G_1(t) + 2 \ln\Bigl(1 + \Bigl|S(t) + \frac{b'(t)}{b(t)}\Bigr|\Bigr)\biggr] < + \infty,}{\sup\limits_{t\ge t_0}\biggl[\il{t_0}{t} \mathfrak{Re} \hskip 2pt S(\tau) d\tau + \ln |b(t)| + \frac{1}{2}\ln G_1(t)\biggr] < + \infty.}  \eqno (3.5)
$$
$$
\left(\sist{\lim\limits_{t\to +\infty}\biggl[\il{t_0}{t} \mathfrak{Re} \hskip 2pt S(\tau) d\tau + \ln |b(t)| - \frac{1}{2}\ln G_1(t) +2 \ln\Bigl(1 + \Bigl|S(t) + \frac{b'(t)}{b(t)}\Bigr|\Bigr)\biggr] = - \infty,}{\lim\limits_{t\to +\infty}\biggl[\il{t_0}{t} \mathfrak{Re} \hskip 2pt S(\tau) d\tau + \ln |b(t)| + \frac{1}{2}\ln G_1(t)\biggr]=- \infty.}\right) ~ \eqno (3.6)
$$
}

Proof. By virtue of Theorem 2.2 it follows from the conditions 9), 10) that Eq. (2.11) is Lyapunov (asymptotically) stable if and only if the inequalities (3.5) (the equalities (3.6)) hold. Then by Theorem 2.5 (statement II) from 8) it follows that the system (1.1) is Lyapunov (asymptotically) stable if and only if the inequalities (3.5) (the equalities (3.6)) are satisfied. The theorem is proved.

By analogy can be proved

{\bf Theorem 3.5.} {\it Let the condition 8) of Theorem 3.4 and the following conditions be satisfied

11) $G_1(t) < 0, \ph t\ge t_0, \ph S(t) + \frac{b'(t)}{b(t)}$ and $G_1(t)$ are continuously differentiable on $[t_0,+\infty)$;

12) $G_1(t)$ is non increasing and for some $\varepsilon > 0$ the function $\frac{G_1'(t)}{|G_1(t)|^{3/2 - \varepsilon}}$ is bounded or

$12_1$) $-G_1(t) \ge\varepsilon > 0$, the function $\frac{G_1'(t)}{G_1(t)}$ is bounded and $\ilp{t_0}\rho_{|G_1|}(\tau) \frac{|G_1'(\tau)|}{|G_1(\tau)|^{3/2}} d \tau < + \infty.$

Then the system (1.1) is Lyapunov (asymptotically) stable if and only if
$$
\sup\limits_{t\ge t_0}\biggl[\il{t_0}{t}\biggl(\mathfrak{Re}\hskip 2pt S(\tau) + 2\sqrt{|G_1(\tau)}\biggr) d \tau + \ln |b(t)| + \phantom{aaaaaaaaaaaaaaaaaaaaaaaaaaaaaaa}
$$
$$
\phantom{aaaaaaaaaaaaaaaaaaaaa}+2\ln \Bigl[1 + \Bigl|S(t) + \frac{b'(t)}{b(t)} + 2\sqrt{|G_1(t)|}\Bigr|\Bigr] - \frac{1}{2} \ln |G_1(t)|\biggr] < + \infty
$$
$$
\biggl(\lim\limits_{t\to + \infty}\biggl[\il{t_0}{t}\biggl(\mathfrak{Re}\hskip 2pt S(\tau) + 2\sqrt{|G_1(\tau)}\biggr) d \tau + \ln |b(t)| + \phantom{aaaaaaaaaaaaaaaaaaaaaaaaaaaa}\biggr.
$$
$$
\biggl.\phantom{aaaaaaaaaaaaaaaaaaaa}+2\ln \Bigl[1 + \Bigl|S(t) + \frac{b'(t)}{b(t)} + 2\sqrt{|G_1(t)|}\Bigr|\Bigr] - \frac{1}{2} \ln |G_1(t)|\biggr] = - \infty.\biggr)
$$
}
$\phantom{aaaaaaaaaaaaaaaaaaaaaaaaaaaaaaaaaaaaaaaaaaaaaaaaaaaaaaaaaaaaaaaaaaaaaaa} \Box$

{\bf Remark 3.2.} {\it On the basis of Corollary 2.1 and Theorem 2.6 one can conclude that the conditions 7) and $7_1)$
of Theorem 3.3 can be replaced by the following simple ones.
$$
-G_k(t) \ge \varepsilon > 0, \ph t\ge t_0, \ph \frac{|G_k'(t)|}{|G_k(t)|} \le \frac{M}{(1 + t - t_0)} \alpha_k, \ph t\ge t_0, \ph \alpha_k > 0,
$$
$$
\ilp{t_0}\frac{d\tau}{\sqrt{|G_k(\tau)|}(1 + \tau - t_0)^{2\alpha_k}} < + \infty, \ph k=1,2.
$$
Similar conclusions are valid with respect to the conditions of Theorem 3.2, Theorem 3.4 and Theorem 3.5.
}

{\bf Remark 3.3.} {\it Let $a_0, \ph b_0, \ph  c_0$ and $d_0$ be real constants. Consider the linear system
$$
\sist{\phi' = a_0 \phi + b_0 \psi,}{\psi' =c_0 \phi + d_0 \psi, \ph t\ge t_0.}
$$
According to the Routh - Hurwitz's  criterion (see [1], pp. 105, 106) this system is asymptotically stable if and only if
$$a_0 + d_0 < 0\phh and \phh  a_0 d_0 - b_0 c_0 > 0.$$
\noindent
Then it is not difficult to verify that (except the trivial cases $G_1(t) = G_2(t) \equiv 0$ and $b(t)= c(t) \equiv 0$) in the two dimensional case the Routh - Hurwitz's  criterion is a consequence of the group of Theorem 3.4 and Theorem 3.5 (in these theorems the restrictions on $c(t)$ are not obligatory [see Remark 2.2]).}

\vskip 20 pt

\centerline{References}

\vskip 20 pt

\noindent
1. L. Y. Adrianoba,  Introduction to the theory of linear systems of differential equations. \linebreak \phantom{a}
S. Peterburg, Publishers of St. Petersburg University, 1992.

\noindent
2. G. A. Grigorian, On the Stability of Systems of Two First - Order Linear Ordinary\linebreak \phantom{a} Differential Equations, Differ. Uravn., 2015, vol. 51, no. 3, pp. 283 - 292.

\noindent
3. G. A. Grigorian. Necessary Conditions and a Test for the Stability of a System of Two\linebreak \phantom{a} Linear Ordinary Differential Equations of the First Order. Differ. Uravn., 2016, \linebreak \phantom{aa} Vol. 52, No. 3, pp. 292 - 300.

\noindent
4. G. A. Grigoryan, “Stability Criterion for Systems of Two First-Order Linear Ordinary \linebreak \phantom{aa}  Differential Equations”, Mat. Zametki, 103:6 (2018),  831–840.

\noindent
5. L. Cezary,   Asymptotic behavior and stability of solutions of ordinary differential \linebreak \phantom{aa}  equations.  Moscow, ''mir'', 1964.

\noindent
6. V. A. Yakubovich, V. M. Starzhinsky,  Linear differential equations with periodic \linebreak \phantom{aa}  coefficients and their applications. Moscow, ''Nauka'', 1972.

\noindent
7.  Ph. Hartman, Ordinary differential Equations.
Second Edition, SIAM, 2002.

\noindent
8. R. Bellman,  Stability theory of differential equations.
 Moscow, Foreign Literature \linebreak \phantom{aa} Publishers, 1954.

\noindent
9. M. V. Fedoriuk. Asymptotic methods for linear ordinary differential equations. \linebreak \phantom{aa} Moskow, ''Nauka'', 1983.

\noindent
10. I. M. Sobol. Study of the asymptotic behaviour of the solutions of the linear \linebreak \phantom{aa} second order differential equations wit the aid of polar coordinates. "Matematicheskij \linebreak \phantom{aa} sbornik", vol. 28 (70), N$^\circ$ 3, 1951, pp. 707 - 714.

\noindent
11. G. A. Grigorian. Boundedness and stability criteria for linear ordinary differential \linebreak \phantom{aa} equations of the second order. "Izvestia vuzov, Matematika, $N^\circ$ 12, 2013, pp. 11 - 18.

\noindent
12.  G. A. Grigorian.  Stability criteria for   second order linear ordinary differential equations. \linebreak \phantom{aa} Saraevo J. Math. In print.

\end{document}